\documentclass[11pt]{amsart}

\usepackage{times,epsf,amssymb,amsmath,rlepsf,hyperref}

\begin{document}

\newtheorem{thm}{Theorem}[section]
\newtheorem{lem}[thm]{Lemma}
\newtheorem{cor}[thm]{Corollary}

\theoremstyle{definition}
\newtheorem{defn}{Definition}[section]

\theoremstyle{remark}
\newtheorem{rmk}{Remark}[section]

\def\square{\hfill${\vcenter{\vbox{\hrule height.4pt \hbox{\vrule
width.4pt height7pt \kern7pt \vrule width.4pt} \hrule height.4pt}}}$}

\newenvironment{pf}{{\it Proof:}\quad}{\square \vskip 12pt}

\title{Non-properly Embedded Minimal Planes in Hyperbolic 3-Space}
\author{Baris Coskunuzer}
\address{Department of Mathematics \\ Koc University \\ Sariyer, Istanbul 34450 Turkey}
\email{bcoskunuzer@ku.edu.tr}
\thanks{The author is partially supported by EU-FP7 Grant IRG-226062, TUBITAK Grant 109T685 and TUBA-GEBIP Award.}

\maketitle

%% User definitions:

\newcommand{\Si}{S^2_{\infty}(\mathbf{H}^3)}
\newcommand{\PI}{\partial_{\infty}}
\newcommand{\SI}{S^{n}_{\infty}(\mathbf{H}^{n+1})}
\newcommand{\BHH}{\mathbf{H}^{n+1}}
\newcommand{\CH}{\mathcal{C}(\Gamma)}
\newcommand{\BH}{\mathbf{H}^3}
\newcommand{\BR}{\mathbf{R}}
\newcommand{\BC}{\mathbf{C}}
\newcommand{\BZ}{\mathbf{Z}}
\newcommand{\T}{\mathcal{T}}
\newcommand{\A}{\mathcal{A}}
\newcommand{\D}{\mathcal{D}}

\begin{abstract}

In this paper, we show that there are non-properly embedded minimal surfaces with finite topology in a simply connected
Riemannian $3$-manifold with nonpositive curvature. We show this result by constructing a non-properly embedded minimal
plane in $\BH$. Hence, this gives a counterexample to Calabi-Yau conjecture for embedded minimal surfaces in negative
curvature case.

\end{abstract}

\section{Introduction}

In the last decade, there have been a great progress in the understanding of the global structure of complete, embedded
minimal surfaces in $\BR^3$. In \cite{CM}, Colding and Minicozzi showed that a complete, embedded minimal surfaces with
finite topology in $\BR^3$ is proper. This result is also known as Calabi-Yau conjecture for embedded minimal surfaces
in the literature. Later, Meeks and Rosenberg generalized this result to complete, embedded minimal surfaces with
positive injectivity radius in $\BR^3$ in \cite{MR1}. A nice survey on Calabi-Yau problem can be found at \cite{Al}.

If we consider the question in more general settings, it is not hard to construct examples of non-properly embedded minimal surfaces with finite
topology in a non-simply connected $3$-manifold, or in a simply connected $3$-manifold with some arbitrary metric. However, there is no known example
of a non-properly embedded minimal surface with finite topology in a simply connected $3$-manifold with nonpositive curvature (See Final Remarks).

On the other hand, the key lemma of \cite{MR1} to prove its main result also applies to all $3$-manifolds with
nonpositive curvature. Hence, the question of whether the generalization of Calabi-Yau conjecture for embedded minimal
surfaces to simply connected $3$-manifolds with nonpositive curvature is true or not became interesting (See Final
Remarks). In this paper, we construct an example of non-properly embedded minimal plane in $\BH$ which shows that the
Calabi-Yau conjecture for embedded minimal surfaces does not generalize to simply connected $3$-manifolds with
nonpositive curvature.\\

\noindent {\bf Theorem 2.1} There exists a non-properly embedded, complete minimal plane in $\BH$.\\

This example is inspired by the heuristic construction of Freedman and He mentioned also in \cite{Ga}. Note that Meeks
and Perez also mentioned such an example in their survey paper \cite{MP}. The idea is as follows: Take a sequence of
round circles $\{C_n\}$ in $\Si$ which limits on the equator circle. Each circle $C_n$ bounds a geodesic plane $P_n$.
By connecting each circle $C_n$ with $C_{n+1}$ by using bridges in $\Si$, we get a nonrectifiable curve $\Gamma$ in
$\Si$ (See Figure 1). Then, we construct a special sequence of minimal disks $\{E_n\}$ with $\partial E_n \to \Gamma$
and show that the limiting minimal plane $\Sigma$ with $\PI \Sigma = \Gamma$ does not stay close to $\Si$ by using
barrier tunnels. Then, we prove that $\Sigma$ is a non-properly embedded minimal plane in $\BH$. Intuitively, one can
imagine $\Sigma$ as the collection of geodesic planes $\{P_n\}$ which are connected via bridges at infinity (See Figure
3).

\subsection{Acknowledgements}

I would like to thank Michael Freedman, William Meeks and David Gabai for very useful conversations, and valuable comments.

\section{The Construction}

First, we need a few definitions. Basic notions and results which will be used in this paper can be found in the survey article \cite{Co2}.

\begin{defn}A {\em least area disk} in a Riemannian manifold $M$ is a compact disk which has the smallest area among the disks in $M$ with the same
boundary. A {\em least area plane} is a complete plane such that any compact subdisk in the plane is a least area disk.
A {\em least area annulus} (compact case) is the annulus which has the least area among the annuli with the same
boundary. A {\em minimal surface} is a surface whose mean curvature vanishes everywhere.
\end{defn}

We will finish the construction in 3 steps. First, we construct a special sequence of compact minimal disks $\{E_n\}$
which will limit on a minimal plane. In second part, we show that the limit of $\{E_n\}$ give us a complete, embedded
minimal plane $\Sigma$ in $\BH$. Finally, we show that $\Sigma$ is indeed non-properly embedded in $\BH$.

\subsection{The Sequence}

\ \\

In this step, we will construct a sequence of minimal disks $\{E_n\}$ in $\BH$ which will give us a complete embedded
minimal plane in $\BH$ as the limit in the next step.

First, consider the half space model for $\BH$. Here $\Si = \BR^2\times \{0\} \cup \{\infty\}$. Now, define a sequence
of round circles $\{C_n\}$ in $\Si$ such that $C_n$ is the round circle in $\BR^2\times \{0\}$ with center at origin
and radius $r_n=1+\frac{1}{n}$. Hence, the sequence limits on the unit circle in $\BR^2\times \{0\}$.

Now, we want to connect each consecutive circle $C_n$ and $C_{n+1}$ with thin bridges in $\Si$. Let's start with $C_1$
and $C_2$. Like in \cite{Ha}, we construct a tunnel in $\BH$ as follows. Let $\eta_1^+$ and $\eta_1^-$ be sufficiently
small round circles in $\BR^2\times\{0\}$ with radius $\delta_1$ where $2\delta_1 < r_1 - r_2$. Let the center of
$\eta_1^+$ be $(\frac{r_1+r_2}{2},\epsilon_1)$ and let the center of $\eta_1^-$ be $(\frac{r_1+r_2}{2},-\epsilon_1)$ in
$\BR^2$ where $\epsilon_1>\delta_1>0$. Let $P_1^\pm$ be the geodesic plane in $\BH$ with $\PI P_1^\pm = \eta_1^\pm$. By
changing $\epsilon_1$ and $\delta_1$ if necessary, we can find sufficiently large disk $D_1^\pm$ in $P_1^\pm$ with
$\partial D_1^\pm = \beta_1^\pm$ such that $\beta_1^+ \cup \beta_1^-$ bounds a least area annulus $A_1$ in $\BH$ like
in \cite{Ha}. Define the annulus $\A_1$ such that $\A_1=(P_1^+-D_1^+)\cup (P_1^- - D_1^-) \cup A_1$. Then, $\A_1$
separates $\BH$ into two components. Let $\mathcal{T}_1$ be the component with $\PI \mathcal{T}_1 =
\Delta_1^+\cup\Delta_1^-$ where $\Delta_1^\pm$ is the disk in $\BR^2$ with boundary $\eta_1^\pm$. We will call
$\mathcal{T}_1$ as a {\em tunnel}. Note that $\partial \T_1 = \A_1$ and $\PI \T_1 = \Delta_1^+ \cup \Delta_1^-$.

Now, we will connect $C_1$ and $C_2$ with a bridge. Let $\mu_1$ be an arc segment in $C_1$ with $\mu_1= C_1 \cap \left
( [1,3]\times (-\frac{\epsilon_1-\delta_1}{2}, \frac{\epsilon_1-\delta_1}{2}) \right )$ and $\mu_2$ be an arc segment
in $C_2$ with $\mu_2= C_2 \cap \left ( [1,3]\times (-\frac{\epsilon_1-\delta_1}{2},
\frac{\epsilon_1-\delta_1}{2})\right )$. Let the end points of $C_n-\mu_n$ be $p^+_n$ and $p^-_n$ where $p^+_n$ be the
endpoint belonging to the upper half space in $\BR^2$. Let $l_1^+$ be the straight line segment in $\BR^2$ between
$p^+_1$ and $p^+_2$. Likewise, let $l_1^-$ be the straight line segment in $\BR^2$ between $p^+_1$ and $p^+_2$. Now,
define $C_1 \#  C_2$ such that $C_1 \#  C_2 = (C_1-\mu_1) \cup (C_2-\mu_2) \cup l^+_1 \cup l^-_1$. Let $\Gamma_1 = C_1$
and $\Gamma_2 = C_1  \#  C_2$. Notice that $\Gamma_2$ separates $\Si$ into two parts, where one of the parts contains
$\Delta_1^+$ and $\Delta_1^-$ which we call the base of the tunnel $\T_1$. In the Poincare ball model for $\BH$, one
can think that the bridge constructed with $l_1^+\cup l^-_1$ goes over the tunnel $\mathcal{T}_1$.

\begin{figure}[t]

\relabelbox  {\epsfxsize=5in

\centerline{\epsfbox{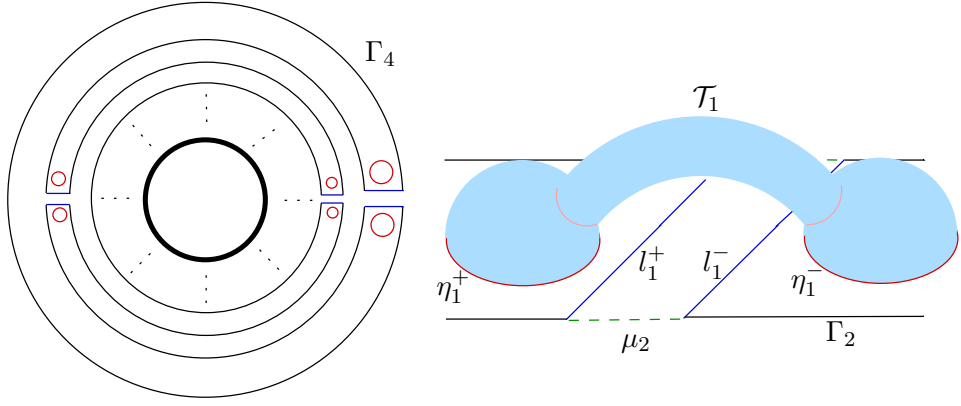}}}

\relabel{1}{}

\relabel{2}{}

\relabel{3}{$\Gamma_4$}

\relabel{4}{}

\relabel{5}{$l_1^+$}

\relabel{6}{$\l_1^-$}

\relabel{7}{$\mu_2$}

\relabel{8}{$\Gamma_2$}

\relabel{9}{$\T_1$}

\relabel{10}{$\eta_1^+$}

\relabel{11}{$\eta_1^-$}

\endrelabelbox

\caption{\label{fig:figure2} \small {In the figure left, $\Gamma_4$ is constructed by using the circles $C_1, C_2, C_3,
C_4$ in $\Si$ by connecting each other with bridges (blue line segments $l_i^\pm$). The red circles represents
$\eta_i^\pm$. In the figure right, the tunnel $\T_1$ is shown.}}

\end{figure}

Next, we will connect $\Gamma_2$ and $C_3$ with a bridge. This time, the bridge will be constructed in the opposite
side, i.e. near the line segment $[-r_3,-r_2]\times \{0\}$ line segment in negative $x$-axis of $\BR^2$. Let $\lambda_1
= \frac{r_2-r_3}{r_1-r_2}$. Now, do the same construction described in the previous paragraphs for
$\epsilon_2=\lambda_1 \cdot \epsilon_1$ and $\delta_2= \lambda_1 \cdot \delta_1$. Here, $\eta_2^\pm$ will be the round
circles in $\BR^2$ with radius $\delta_2$ and center $(-\frac{r_2+r_3}{2}, \pm \epsilon_2)$. Let $\phi_2$ be the
isometry of $\BH$ with $\phi_1(\eta_1^\pm)=\eta_2^\pm$ and $\phi_1((\frac{r_1+r_2}{2},0))=(\frac{r_2+r_3}{2},0)$.
Notice that $\phi_2$ can be obtained by composing the parabolic isometry which translates the point
$(\frac{r_1+r_2}{2},0)$  to the point $(\frac{r_2+r_3}{2},0)$ fixing $\{\infty\}$, and the hyperbolic isometry with
dilation constant $\lambda_1$. Then, define the second tunnel $\T_2$ as the isometric image of $\T_1$, i.e.
$\T_2=\phi_2(\T_1)$. Define annuli $A_2$ and $\A_2$ accordingly. The second bridge can be defined similarly. Let
$\tau_2= C_2 \cap \left ( [-3,-1]\times (-\frac{\epsilon_2-\delta_2}{2}, \frac{\epsilon_2-\delta_2}{2})\right )$ and
$\tau_3= C_3 \cap \left ( [-3,-1]\times (-\frac{\epsilon_2-\delta_2}{2}, \frac{\epsilon_2-\delta_2}{2})\right )$. For
$n \geq 2$, let the end points of $C_n-\tau_n$ be $q^+_n$ and $q^-_n$ where $q^+_n$ be the endpoint belonging to the
upper half space in $\BR^2$. Let $l_2^+$ be the straight line segment in $\BR^2$ between $q^+_2$ and $q^+_3$.
Similarly, $l_2^-$ be the straight line segment in $\BR^2$ between $q^-_2$ and $q^-_3$. Like before define $\Gamma_2 \#
C_3$ such that $\Gamma_2 \#  C_3= (\Gamma_2- \tau_2) \cup (C_3 - \tau_3) \cup l_2^+ \cup l_2^-$. Let $\Gamma_3=
\Gamma_2  \#  C_3$.

Hence, if we continue the process, we can define the tunnels and bridges for any $n$ as follows. Let $o_n$ be the
largest odd number which is smaller than or equal to $n$. Let $e_n$ be the largest even number which is smaller than or
equal to $n$. In other words, if $n$ is odd, $o_n=n$, and if $n$ is even $o_n= n-1$. Similarly, if $n$ is even, $e_n=n$
and if $n$ is odd, $e_n=n-1$.

Now, let's define the basic components in the construction of bridges and tunnels in general terms. Let $\lambda_n=
\frac{r_{n+1}-r_{n+2}}{r_1-r_2}$. Then define $\epsilon_n= \lambda_n \cdot \epsilon_1$ and $\delta_n = \lambda_n \cdot
\delta_1$. $\mu_n$ is an arc segment in the circle $C_n$ such that $\mu_n= C_n \cap \left ( [1,3] \times
(-\frac{\epsilon_{o_n}-\delta_{o_n}}{2}, \frac{\epsilon_{o_n}-\delta_{o_n}}{2})\right )$.  Similarly, $\tau_n$ is an
arc segment in the circle $C_n$ such that $\tau_n = C_n \cap \left ( [-3,-1] \times (-\frac{\epsilon_{e_n} -
\delta_{e_n}}{2}, \frac{\epsilon_{e_n} - \delta_{e_n}}{2})\right )$. Then when $n$ is odd, $l_n^\pm$ would be the
straight line segments between the points $p_n^\pm$ and $p_{n+1}^\pm$, and when $n$ is even, $l_n^\pm$ would be the
straight line segments between the points $q_n^\pm$ and $q_{n+1}^\pm$. Hence, the arc segments $\mu_n$, and odd indexed
line segments $l_n^\pm$ live in the right side ($x>0$) of $\BR^2$, and the arc segments $\tau_n$ and even indexed line
segments $l_n^\pm$ live in the left side ($x<0$) of $\BR^2$.

Define $\Gamma_n$ inductively such that $\Gamma_{n+1} = \Gamma_n  \#  C_{n+1}$. Here $\Gamma_n  \#  C_{n+1}$ can be
defined as follows. When $n$ is odd, $\Gamma_n  \#  C_{n+1} = (\Gamma_n-\mu_n) \cup (C_{n+1} - \mu_{n+1}) \cup l_n^+
\cup l_n^-$. When $n$ is even, $\Gamma_n  \#  C_{n+1} = (\Gamma_n-\tau_n) \cup (C_{n+1} - \tau_{n+1}) \cup l_n^+ \cup
l_n^-$. By iterating the procedure, we will get a non-rectifiable connected curve $\Gamma$ (or $\Gamma_\infty$) in
$\Si$ which has infinite length.

We define the tunnels as follows. When $n$ is odd, the round circles $\eta_n^\pm$ would be the round circles of radius
$\delta_n$ with center $(\frac{r_n+r_{n+1}}{2},\pm \epsilon_n)$. Similarly, when $n$ is even, the round circles
$\eta_n^\pm$ would be the round circles of radius $\delta_n$ with center $(-\frac{r_n+r_{n+1}}{2},\pm \epsilon_n)$.
Then $P_n^\pm$ would be the geodesic planes in $\BH$ with $\PI P_n^\pm = \eta_n^\pm$. Define the least area annulus
$A_n$ such that the boundary curves are in the geodesic planes $P_n^+$ and $P_n^-$. Hence, the tunnel $\mathcal{T}_n$
can be defined accordingly by using $P_n^\pm$ and $A_n$. Notice that when $n$ is odd, the tunnels $\mathcal{T}_n$ (and
hence $\eta_n^\pm$, $P_n^\pm$, $A_n$) are in the right side ($x>0$) in $\BR^2$, and when $n$ is even, the tunnels
$\mathcal{T}_n$ (and hence $\eta_n^\pm$, $P_n^\pm$, $A_n$) are in the left side ($x<0$) in $\BR^2$. In other words, odd
indexed tunnels and bridges are in the right side of $\BR^2$ and even indexed tunnels and bridges are in the left side
of $\BR^2$.

Now, we construct the sequence of embedded minimal disks $\{E_n\}$ in $\BH$. Let $X_1= \BH- int(\mathcal{T}_1)$. Then
$\partial X_1 = \A_1$ and $\PI X_1 = \Si - (\Delta_1^+\cup \Delta_1^-)$. In general, let $X_n = \BH - \bigcup_{i=1}^n
int(\T_i)$. Then, $\partial X_n = \bigcup_{i=1}^n \A_n$ and $\PI X_n = \Si- \bigcup_{i=1}^n int(\Delta_i^+\cup
\Delta_i^-)$. Define $X$ such that $X=\bigcap_{i=1}^\infty X_n = \BH - \bigcup_{i=1}^\infty int(\T_i)$

Let $\upsilon$ be the circle of radius $3$ with center at origin in $\BR^2$. Let $\mathcal{P}$ be the geodesic plane in
$\BH$ with $\PI \mathcal{P} = \upsilon$. Then $\mathcal{P}$ separates $\BH$ into two components. Let $\Omega$ be the
component whose asymptotic boundary contains the origin. Also, define the horoball $\mathcal{H}_i$ such that
$\mathcal{H}_i = \{ \ (x,y,z) \in \BH \ | \ z\geq \frac{1}{i} \ \}$. Let $S_i$ be the horosphere with $S_i = \partial
H_i$, i.e. $S_i = \{ \ (x,y,\frac{1}{i}) \in \BH \ \}$. Then define the domain $\Omega_i$ such that $\Omega_i = \Omega
\cap \mathcal{H}_i \cap X$.  Notice that $\Omega_i$ is a compact mean convex domain since the horosphere $S_i$ has mean
curvature $1$, and the geodesic plane $\mathcal{P}$ and the least area annulus $\A_n$ have mean curvature $0$.

Let $p$ be the point $(0,0,1)$ in the upper half space model of $\BH$. Notice that $\{p\} = \Sigma_1\cap l$ where $\Sigma_1$ is the geodesic plane in $\BH$ with $\PI \Sigma_1 = \Gamma_1$, and $l$ is the vertical geodesic $l= \{ \ (0,0,t) \in \BH \ |
\ t>0 \ \}$. Let $\mathcal{C}_n$ be the geodesic cone over $\Gamma_n$ in $\Si$ with
cone point $p$. In other words, if $\gamma_{pq}$ is the geodesic ray in $\BH$ starting from $p$ and limiting on $q\in
\Si$, then $\mathcal{C}_n = \bigcup_{q\in\Gamma_n} \gamma_{pq}$. Let $\alpha_n^i$ be a simple closed curve in $\BH$
defined as the intersection of the horosphere $S_i$ and the geodesic cone $\mathcal{C}_n$, i.e. $\alpha_n^i = S_i\cap
\mathcal{C}_n$. Notice that for any $n$, there exists a number $c_n>0$ such that if $i\geq c_n$ then $\alpha_n^i \cap
\T_n = \emptyset$. Now, we need a lemma due to Meeks-Yau in order to continue to the construction.

\begin{lem} \cite{MY}
Let $\Omega$ be a compact, mean convex $3$-manifold, and $\alpha\subset\partial \Omega$ be a nullhomotopic simple
closed curve. Then, there exists an embedded least area disk $ D\subset M$ with $\partial  D = \Gamma$.
\end{lem}

By construction, for any $i$, $\Omega_i$ is a compact, mean convex $3$-manifold. Also, for any $n$, $\alpha_n^{c_n}$ is
a simple closed curve in $\partial \Omega_{c_n}$. Hence, by Lemma 2.1, there exists an embedded least area disk $E_n$
in $\Omega_{c_n}$ such that $\partial E_n = \alpha_n^{c_n}$ (say $\alpha_n$ for short). Note that being least area in
$\Omega_{c_n}$ does not imply that $E_n$ is also least area in $\BH$. Even though $E_n$ may not be not a least area
disk in $\BH$, we claim that it is indeed least area in $X= \BH - \bigcup_{i=1}^\infty int(\T_i)$.

\begin{lem} For any $n$, $E_n$ is a least area disk in $X$.
\end{lem}

\begin{pf} We know that $E_n$ is a least area disk in $\Omega_{c_n}$ by construction. The only case where $E_n$ is not a
least area disk in $X$ is the existence of another disk $\mathcal{D}$ in $X$ with $\partial \mathcal{D}= \partial E_n =
\alpha_n$ and $|\mathcal{D}| < |E_n|$ where $|.|$ represents the area. Notice that the horoball $\mathcal{H}_i$
and the half space $\Omega$ are convex subsets of $\BH$. Hence, for any $i$, $\widehat{\Omega}_i = \mathcal{H}_i \cap
\Omega$ is a convex subset of $\BH$. Then if $\pi_i:\BH\to \widehat{\Omega}_i$ is the nearest point projection, since
$\widehat{\Omega}_i$ is convex, $\pi_i$ would be distance reducing. Let $\widehat{\D}=\pi_{c_n} (\mathcal{D})$ (say $\pi$ for $\pi_{c_n}$ for short). Then,
$|\widehat{\D}| \leq |\mathcal{D}|$ and $\widehat{\D}$ is a disk in $\widehat{\Omega}_{c_n}$ with $|\widehat{\D}|\leq
|\D| < |E_n|$.

The outline of the proof is as follows: If we can show that $\widehat{\D}$ is disjoint from any tunnel $\T_i$, then we
get a contradiction, as $\widehat{\D}$ would a be a disk in $\Omega_{c_n}$ with $\partial \widehat{\D} =
\alpha_n$, and $|\widehat{\D}| < |E_n|$. Hence, if $\widehat{\D}$ intersects some tunnels $\T_i$, we will do a
surgery on $\widehat{\D}$ by removing the subdisks where $\widehat{\D}$ intersect the tunnels, and get a new disk
$\widehat{\D}'$ in $\Omega_{c_n}$ with $\partial \widehat{\D}' = \alpha_n$ with $|\widehat{\D}'| < |E_n|$. Since
$E_n$ is the least area disk in $\Omega_{c_n}$ with boundary $\alpha_n$, this will give us a contradiction.

Notice that $\pi(x)= x$ for any $x \in \D \cap int(\Omega_{c_n})$, and $\pi (\D - int(\Omega_{c_n}))
\subset S_{c_n}$. Let $T_i\cap S_{c_n} = O_i^+ \cup O_i^-$. Hence, if we can show that $\D \cap O^\pm_i = \emptyset$,
we are done.

First, let $\psi: D^2 \to X$ be a parametrization of $\D$, i.e. $\psi(D^2) =\D$. Then let $\varphi: D^2 \to
\Omega_{c_n}$ be the parametrization of $\widehat{\D}$ with $\varphi = \pi \circ \psi$, i.e.
$\varphi(D^2)=\widehat{\D}$. If $\widehat{\D}\cap O_i^+ \neq \emptyset$, there are two cases: $O_i^+ \nsubseteq
\widehat{\D}$ or $O^+_i \subset \widehat{\D}$. If $\widehat{\D}\cap O_i^+ \neq \emptyset$ and $O_i^+ \nsubseteq
\widehat{\D}$, then let $V_i^+ = \varphi^{-1}(O_i^+ \cap \widehat{\D})$. Since $\T_i$ is separating, $\partial V_i^+$
would be a collection of circles. Since $O_i^+ \nsubseteq \widehat{\D}$, $\varphi (\partial V_i^+)$ is a collection of
nonessential circles in $\partial \T_i$. Hence, we can push off $O^+_i\cap \widehat{\D}$, and get another disk
$\widehat{\D}'$ with less area.

Now, if $O^+_i \subset \D$, then we claim that $\varphi^{-1}(O^+_i)$ consists of even number of disks after removing
the nonessential circles like in the previous paragraph. Let $y$ be a point in $O_i^+$. $\pi^{-1} (y)$ would be a
geodesic ray $\rho$ starting at $y$ and orthogonal to $S_i$. Also let $\tau$ be an infinite ray in $\T_i$ starting at
$y$ and outside of $\mathcal{H}_{c_n}$. We know that $\D \cap \tau = \emptyset$ as $\D \subset X$.

Let $\Delta$ be a homotopy between the infinite rays $\tau$ and $\rho$. In other words, $\Delta: [0,1] \times [0,1) \to
\BH$ a continuous map such that $\Delta(\{0\}\times [0,1)) = \tau$, $\Delta(\{1\}\times [0,1)) = \rho$ and
$\Delta([0,1] \times\{0\}) = y$ and $\Delta|{[0,1] \times (0,1)}$ is an embedding. We abuse the notation by using
$\Delta$ for its image. We can also assume $\Delta$ is transverse to $\D$. Then $\Delta \cap \D$ would be some
collection of circles and some paths $\{s_j\}$ whose endpoints $\{q_j^+,q_j^-\}$ are in $\rho$ as $\tau \cap \D =
\emptyset$. Each $p_j^\pm=\psi^{-1}(q_j^\pm) \in D^2$ are in different subdisks $W_j^\pm$ in $D^2$ where $\{W_j^\pm\}
\subset\varphi^{-1}(O^+_i) $. As $\varphi = \pi \circ \psi$, there is a one to one correspondence between the points
$p_j^\pm$ which are in $\varphi^{-1}(y)$ and the disks $W_j^\pm$, i.e. $p_j^\pm \in W_j^\pm$. Recall that
$\varphi(p_j^\pm) = y$ and $\varphi (W_j^\pm) = O^+_i$.

Now, we will modify $\widehat{\D}$ by pushing it off from $O_i^+$ and get a new disk $\widehat{\D}'$ in $\Omega_{c_n}$
with less area. For a fixed $j_0$, let $B_{j_0}$ be a subdisk in the interior of $\D$ which contains
$\psi(W_{j_0}^+)\cup \psi(W_{j_0}^-)$, and no other $\psi(W_j^\pm)$ for $j\neq j_0$. Consider the disk
$\widehat{B}_{j_0}=\pi(B_{j_0})$ in $S_{c_n}$. By construction, $O^+_i=\pi(W_j^\pm)$ is in $\widehat{B}_{j_0}$. Now,
$\partial \widehat{B}_{j_0}$ consists of two parts, say $\partial \widehat{B}_{j_0}= \omega_1\cup\omega_2$ where
$\omega_1=\pi(\partial B_{j_0})$, and $\omega_2 = \partial \widehat{B}_{j_0} - \omega_1$. By construction, $\omega_2
\neq \emptyset$. Notice that $\omega_2$ consists of points where the geodesic ray starting from those points orthogonal
to $S_{c_n}$ are tangent to $\D$. Then by using a version of Sard's theorem, we can find a path $\gamma$ connecting
$\omega_2$ and $\partial O_i^+$ such that $\pi^{-1}(x)$ has even number of preimages (assume exactly $2$ for
simplicity) with $\pi^{-1}(x)=\{z_x^+,z_x^-\}$ for any point $x\in \gamma$ except for a finitely many points $\{a_1,
a_2,.., a_k\}$. Hence, $\varphi^{-1}(\partial O_i^+ \cup \gamma)$ contains a circle $\beta$ which bounds a (singular)
disk $F$ in $D^2$ such that $\psi(F) \supseteq W_{j_0}^\pm$. Notice that $\beta = \beta^+ \cup \beta^- \cup
\{\varphi^{-1}(\{a_1, a_2,.., a_k\})\}$ where $\beta^+$ is the arc segment consists of positive preimages $\{z_x^+\}$,
and $\beta^-$ is defined similarly.

Now, consider the disk defined as $\widehat{\D}' = \varphi (D^2 - int(F))=\widehat{D}- int(O_i^+)$. Note also that
$\varphi (D^2 - F)=\widehat{D}- (O_i^+\cup\gamma)$. Since $\varphi(z_x^+)=\varphi(z^-_x)$, we get a continuous map
$\varphi':D^2\to \Omega_{c_n}$ such that $\widehat{\D}' = \varphi'(D^2)$.  This is because when you remove a subdisk
$F$ from $D^2$, and identify two connected subarcs $\beta^+$ and $\beta^-$ in $\beta=\partial F = \beta^+ \cup
\beta^-$, then the new topological object would be a disk again. Since  the collection of points $\{a_1, a_2,.., a_k\}$
has only one preimage, when we do this surgery in the domain disk, it becomes a topological space $\Delta$ which is a
disk with $k$ pair of points $\{b_1^+, b_1^-,..., b_k^+, b_k^-\}$ identified. However, $\Delta$ can be seen as the
image of another map $F:D^2\to\Delta$ with $F(b_i^\pm)=a_i$. Hence, the surgered object $\widehat{\D}'$ can be seen as
an image of a disk with $\partial \widehat{\D}' = \alpha_n$ and $\widehat{\D}' \cap \T_i =\emptyset$. By construction,
there are only finitely many $i>0$ with $\widehat{\Omega}_{c_n} \cap \T_i \neq \emptyset$, i.e. there exists $K_n>0$
such that for any $i>K_n$, $\widehat{\Omega}_{c_n} \cap \T_i = \emptyset$. Hence, if we do this surgery for any $i\leq
K_n$, we get a disk $\widehat{\D}'$ in $\Omega_{c_n}$ with $\partial \widehat{\D}' = \alpha_n$ with $|\widehat{\D}'| <
|E_n|$. However, this is a contradiction since $E_n$ is the least area disk in $\Omega_{c_n}$ with boundary $\alpha_n$.
\end{pf}

\subsection{The Limit}

\ \\

In the previous section, we constructed a sequence of least area disks $\{E_n\}$ in $X$ with $\partial E_n = \alpha_n
\to \Gamma$ where $\Gamma$ is the non-rectifiable curve in $\PI X \subset \Si$ constructed in the previous part. In
this section, we show that the sequence of least area disks $\{E_n\}$ has a subsequence limiting on an embedded least
area plane $\Sigma$ in $X$ with $\PI \Sigma = \Gamma$ by using the techniques of Gabai \cite{Ga}. Since $\Sigma$ is an
embedded least area plane in $X$, it will be an embedded minimal plane in $\BH$. In the next section, we will show that
$\Sigma$ is also non-properly embedded in $\BH$, and prove the main result of the paper.

First, we need a definition which we use in the following part. For details of the results and notions in this section,
see Section 3 in \cite{Ga}.

\begin{defn} The sequence $\{D_{i}\}$ of smooth embedded disks in a Riemannian manifold $X$ \textit{converges} to
the lamination $\sigma$ if

\begin{itemize}
\item $\sigma = \{ \ x=\lim x_i \ | \ \ x_i \in D_i , \{x_i\} \mbox{ is a convergent sequence in } X \
\}$\\

\item $\sigma = \{ \ x=\lim x_{n_i} \ | \ \ x_i \in D_i , \{x_i\} \mbox{ has a convergent subsequence }
\{x_{n_i}\} \mbox{ in } X \ \}$\\

\item For any $x \in \sigma$, there exists a sequence $\{x_i\}$ with $x_i \in D_i$ and $\lim x_i = x$
such that there exist embeddings $f_{i}: D^2 \to D_i$ which converge in the $C^{\infty }$-topology to a smooth
embedding $f:D^2 \to L_{x}$, where $x_i \in f_{i}(Int(D^2))$, and $L_{x}$ is the leaf of $\sigma $ through $x$, and
$x\in f(Int(D^2))$.
\end{itemize}

We call such a lamination $\sigma$ a \textit{$D^2$-limit lamination} \cite{Ga}.

\end{defn}

In other words, $\{D_i\}$ is a sequence of smooth embedded disks such that the set of the limits of all $\{x_i\}$ with
$x_i\in D_i$ and the set of the limits of the subsequences are the same. This is a very strong and essential condition
on $\{D_i\}$ in order to limit on a collection of \textit{pairwise disjoint} embedded surfaces. Otherwise, one
might simply take a sequence such that $D_{2i+1} = \Sigma_1$ and $D_{2i} = \Sigma_2$ where $\Sigma_1$ and $\Sigma_2$
are intersecting disks. Then, without the first condition ($\sigma$ being just the union of limit points), $\sigma=
\Sigma_1 \cup \Sigma_2$ in this case, which is not a collection of pairwise disjoint embedded surfaces. However, the
first condition forces $\sigma$ to be either $\Sigma_1$ or $\Sigma_2$, not the union of them. By similar reasons, this
condition is also important to make sure the embeddedness of the disks in the collection $\sigma$.

Now, we state the following lemma (\cite{Ga}, Lemma 3.3) which is essential for the following part.

\begin{lem} If $\{E_n\}$ is a sequence of embedded least area disks in $\BH$, where $\partial E_n \to \infty$, then
after passing to a subsequence $\{E_{n_j}\}$ converges to a (possibly empty) $D^2$-limit lamination $\sigma$ by least
area planes in $\BH$.
\end{lem}

\begin{figure}[b]
\begin{center}
$\begin{array}{c@{\hspace{1in}}c} \multicolumn{1}{l}{\mbox{\bf }} &
    \multicolumn{1}{l}{\mbox{\bf }} \\ [-0.53cm]

\relabelbox  {\epsfxsize=2in

\epsffile{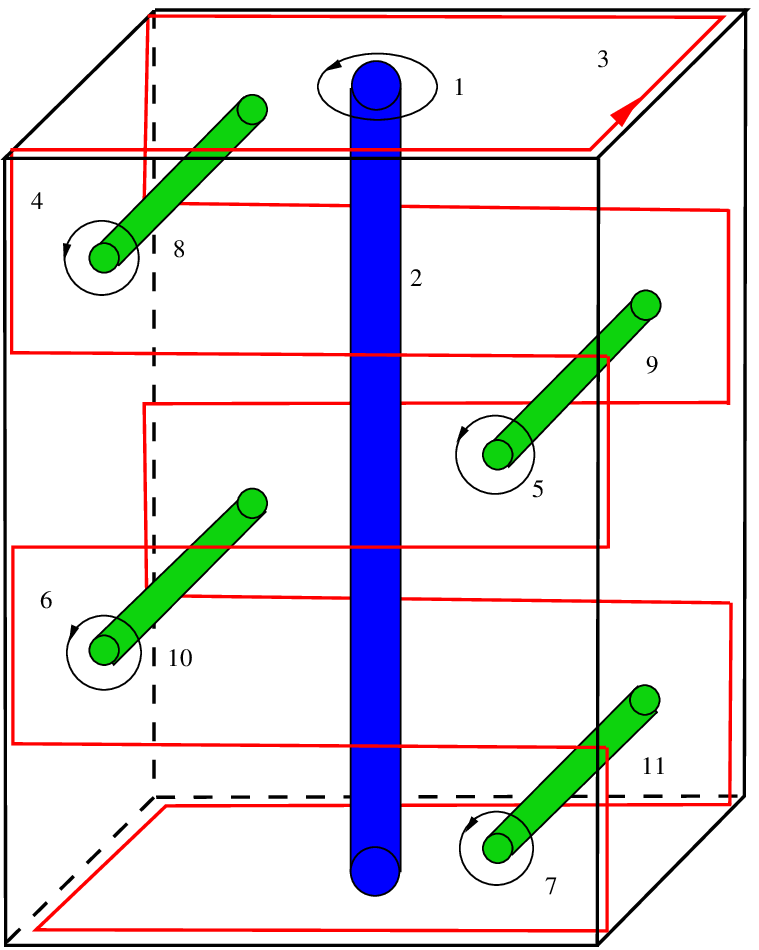}}

\relabel{1}{$\delta$}

\relabel{2}{$\beta$}

\relabel{3}{$\alpha_4$}

\relabel{4}{$\tau_1$}

\relabel{5}{$\tau_2$}

\relabel{6}{$\tau_3$}

\relabel{7}{$\tau_4$}

\relabel{8}{$\T_1$}

\relabel{9}{$\T_2$}

\relabel{10}{$\T_3$}

\relabel{11}{$\T_4$}

\endrelabelbox

& \relabelbox  {\epsfxsize=2in

\epsffile{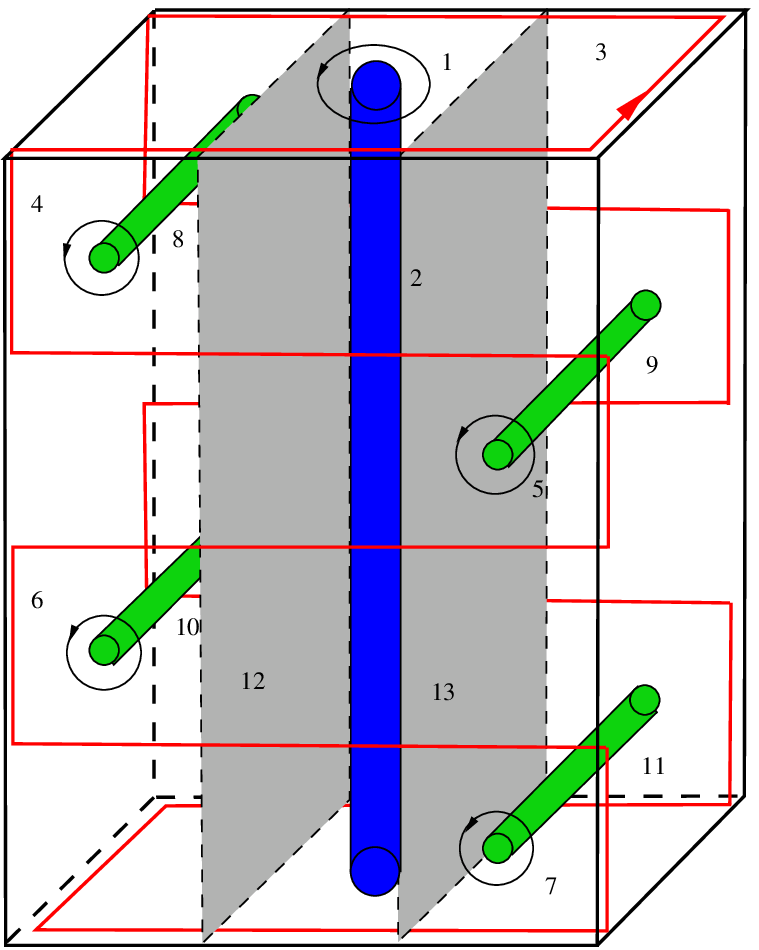}}

\relabel{1}{$\delta$}

\relabel{2}{}

\relabel{3}{$\alpha_4$}

\relabel{4}{$\tau_1$}

\relabel{5}{$\tau_2$}

\relabel{6}{$\tau_3$}

\relabel{7}{$\tau_4$}

\relabel{8}{}

\relabel{9}{}

\relabel{10}{}

\relabel{11}{}

\relabel{12}{$D_1$}

\relabel{13}{$D_2$}

\endrelabelbox \\ [0.4cm]
\mbox{\bf $Z_4$} & \mbox{Collapsing Disks}
\end{array}$
\end{center}
\caption{The box corresponds to $Z_4$. Green tubes are the tunnels, and blue tube is the curve $\beta$. The red curve
is $\alpha_4$. In the right, the grey rectangles are the collapsing disks. } \label{figtest-fig}
\end{figure}

By using this lemma, we will get a least area plane in $X$ which is constructed as a limit of the sequence of least
area disks $\{E_n\}$ in $X$.

\begin{lem} The sequence of least area disks $\{E_n\}$ in $X$ constructed in the previous section has a subsequence $\{E_{n_j}\}$
converges to a nonempty $D^2$-limit lamination $\sigma$ by least area planes in $X$.
\end{lem}

\begin{pf}
Even though Lemma 2.3 is stated for $\BH$, since its proof is a local construction, it also applies to our case where
the ambient manifold is $X \subset \BH$. Hence, all we need to show is that the lamination we get in the limit is
nonempty.

To show that, we construct a sequence of points $\{x_n\}$ with $x_n\in E_n$ such that $\{x_n\}$ has a convergent
subsequence $\{x_{n_j}\}$ with $x_{n_j}\to p$ for $p\in X$. In the upper half space model for $\BH$, let $\beta$ be the
geodesic segment in $\BH$ starting from the point $(0,0, \frac{1}{2})$ and ending at the point $(0,0, 3)$ (See Figure
3). Let $\gamma_1$ be the round circle in $\BR^2\times\{0\}\subset\Si$ with center $(0,0,0)$ and radius $\frac{1}{2}$,
and $\gamma_2$ be the round circle with center $(0,0,0)$ radius $3$. Let $\mathcal{A}$ be the annulus in
$\BR^2\times\{0\}$ bounded by $\gamma_1$ and $\gamma_2$. Then for any $n$, $\Gamma_n$ would be in $\mathcal{A}$. If
$P_i$ is the geodesic plane with $\PI P_i=\gamma_i$, then let $\BH-P_i=\Omega^+_i\cup\Omega^-_i$ where
$(0,0,0)\in\Omega^-_i$. Let $Y= \Omega^+_1\cap\Omega_2^-$. Then, $\PI Y = \mathcal{A}$. Since $Y$ is convex, then by
convex hull property \cite{Co2}, for any $n$, $E_n$ would be in $Y$.

We claim that $E_n\cap \beta \neq \emptyset$ for any $n$. Then by taking $x_n\in \beta \cap E_n$, we can construct the
desired sequence, and finish the proof. We claim that $\partial E_n = \alpha_n$ links $\beta$ in $X$. Recall that for
any $n$, $E_n\subset Y$, and $\partial Y = P_1\cup P_2$. Also, the endpoints of $\beta$ belong to $P_1$ and $P_2$, i.e.
$(0,0,\frac{1}{2})\in P_1$ and $(0,0,3)\in P_2$. Let $Y_n= Y- \bigcup_{i=1}^n \T_i$. Clearly, $Y$ is topologically a
$3$-ball, and $Y_n$ is a genus $n$ handlebody (a $3$-ball with $n$ $1$-handles attached). For any $n$, $\alpha_n$ is a
trivial loop in $\pi_1(Y_n)$ by construction. If we realize $Y_n$ topologically as a $3$-ball with $n$ $1$-handles
attached, and each $1$-handle corresponds to a tunnel $\T_i$, then $\pi_1(Y_n)$ would be free product of $n$ copies of
$\BZ$ i.e. $\pi_1(Y_n)=*_{i=1}^n \BZ$. Hence, $\pi_1(Y_n) = <\tau_1,\tau_2,...,\tau_n>$ where $\tau_i$ is the loop
which corresponds to an essential simple closed curve in the annulus $\partial \T_i$. Again, by construction,
$\alpha_n$ is a trivial loop in $Y_n$.

Let $Z_n=Y_n-\beta$ be the topologically genus $n+1$ handlebody which is a $3$-ball with $n+1$ $1$-handles (See Figure
2). Then, $\pi_1(Z_n)= <\delta, \tau_1,\tau_2, ... , \tau_n>$ where $\delta$ is the generator coming from $\beta$. i.e.
$\delta$ corresponds to the essential loop of the annulus $\partial N_\epsilon(\beta)\cap int(Y)$. Even though,
$\alpha_n$ is trivial in $Y_n$, it is not trivial in $Z_n$ as
$\alpha_n=\delta.\tau_1.\delta^{-1}.\tau_1^{-1}.\tau_2.\delta^{-1}.\tau_2^{-1}....\tau_n.\delta^{-1}.\tau_n^{-1}$ which
is not a trivial element in $\pi_1(Z_n)$. To see this, one might collapse the disks as in Figure 2-right, and divide
$\alpha_n$ to simpler components to write it down explicitly in terms of the generators of $\pi_1(Z_n)$. Hence,
$\alpha_n$ is trivial loop in $Y_n$, but it is not trivial in $Z_n$. This implies that any disk bounding $\alpha_n$ in
$Y_n$ must intersect $\beta$. Hence, $\beta\cap E_n \neq \emptyset$ for any $n$.

Let $x_n$ be a point in $\beta \cap E_n$ for any $n$. The sequence $\{x_n\}$ is a subset of compact geodesic segment
$\beta$. Hence, there is a subsequence $\{x_{n_j}\}$ with $x_{n_j}\to p$ where $p$ is a point in $X$. Now, replace the
sequence $\{E_n\}$ with the subsequence $\{E_{n_j}\}$. Then by applying Lemma 2.3 to the new sequence $\{E_n\}$, we get
a subsequence $\{E_{n_k}\}$ which limits on a nonempty $D^2$-limit lamination $\sigma$ by least area planes in $X$. The
proof follows.

\end{pf}

By the proof of Lemma 2.3, for any leaf $L$ in $\sigma$, for any subdisk $D$ in $L$, we can find sufficiently close
disk $D_n$ in the disk $D_n$. Since $\partial E_n =\alpha_n \to \Gamma$, we can find a least area plane $\Sigma$ in
$\sigma$ with $\PI \Sigma = \Gamma$. Note that again by construction $\PI \sigma = \overline{\Gamma}$ where
$\overline{\Gamma}$ is the closure of $\Gamma$ in $\Si$. By construction of $\Gamma$, $\overline{\Gamma}$ would be
$\Gamma\cup \gamma$ where $\gamma$ is the unit circle in $\BR^2\times\{0\} \subset \Si$ with center $(0,0,0)$. Also, by
varying the transverse geodesic segment $\beta$, it is not hard to show that the geodesic plane $P$ with $\PI P =
\gamma$ is another leaf of the lamination $\sigma$.

\subsection{Non-properly Embeddedness}

\ \\

In this section, we will show that the least area plane $\Sigma$ in $X$ constructed in the previous section is not
properly embedded (See Figure 3). Hence, this will show that $\Sigma$ is a non-properly embedded minimal plane in
$\BH$, and the main result of the paper follows.

\begin{thm}There exists a non-properly embedded, complete minimal plane in $\BH$.
\end{thm}

\begin{pf} We claim that the least area plane $\Sigma$ in $X$ constructed in previous section is not
properly embedded. Since $\Sigma$ is a least area plane in $X$, it is automatically a minimal plane in $\BH$. Hence, if
we show that $\Sigma$ is not properly embedded in $\BH$, we are done.

Assume that $\Sigma$ is properly embedded. Let $\beta$ be as in the proof of Lemma 2.4 (See Figure 3), i.e. In the
upper half space model for $\BH$, $\beta$ is the geodesic segment in $\BH$ starting from the point $(0,0, \frac{1}{2})$
and ending at the point $(0,0, 3)$. If $\beta$ is not transverse to $\Sigma$, modify $\beta$ slightly at non-transverse
points to make it transverse to $\Sigma$. Now, as $\Sigma$ is properly embedded, $\Sigma\cap\beta$ is compact. Let
$\gamma$ be the unit circle in $\BR^2\times\{0\} \subset \Si$ with center $(0,0,0)$, and $P$ be the geodesic plane in
$\BH$ with $\PI P = \gamma$. As mentioned in the previous section $P$ is also a least area plane in the lamination
$\sigma$, and hence $\Sigma\cap P= \emptyset$. Since $\Sigma\cap\beta$ is compact, this implies
$\delta=\inf_{p\in\Sigma\cap\beta} \{\ p_z \in \BR \ | \ p=(p_x,p_y,p_z) \ \}>1$.

\begin{figure}[b]
\begin{center}
$\begin{array}{c@{\hspace{1in}}c} \multicolumn{1}{l}{\mbox{\bf }} &
    \multicolumn{1}{l}{\mbox{\bf }} \\ [-0.53cm]

\relabelbox  {\epsfxsize=2in

\epsffile{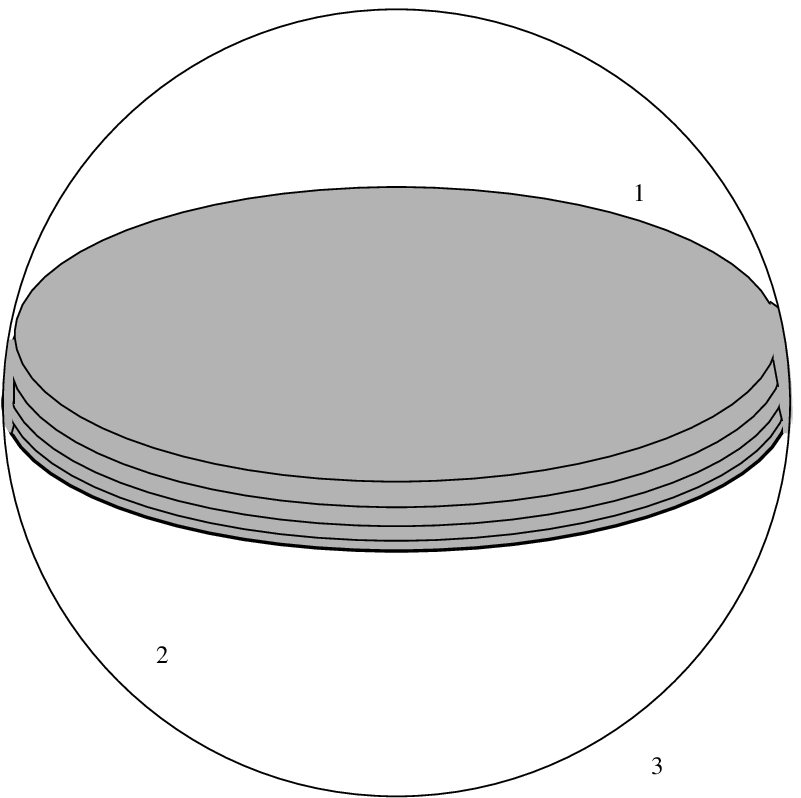}}

\relabel{1}{$\Sigma$}

\relabel{2}{$\BH$}

\relabel{3}{$\Si$}

\endrelabelbox

& \relabelbox  {\epsfxsize=2in

\epsffile{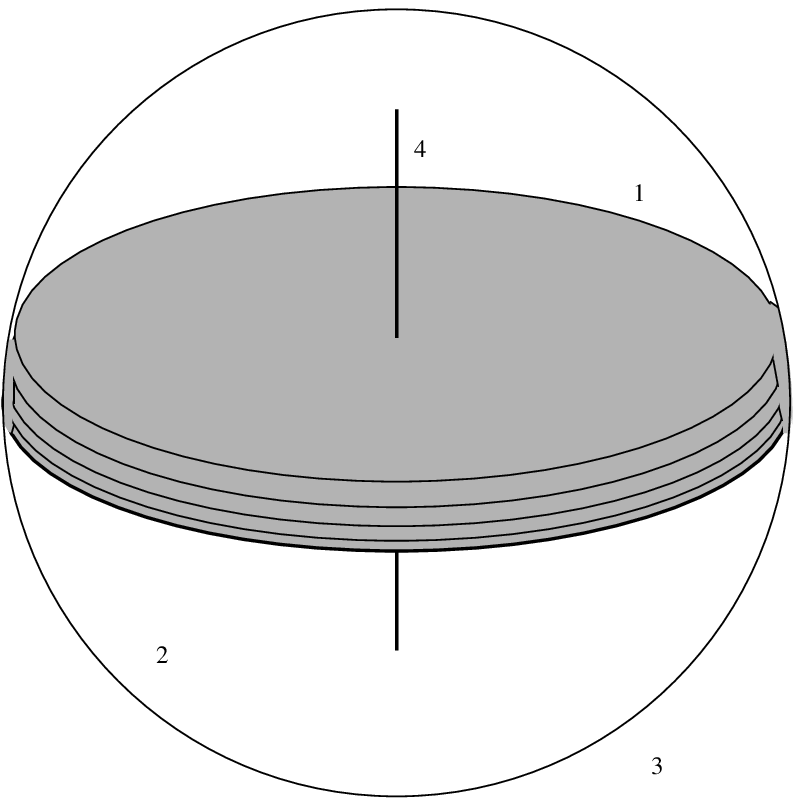}}

\relabel{1}{$\Sigma$}

\relabel{2}{$\BH$}

\relabel{3}{$\Si$}

\relabel{4}{$\beta$}

\endrelabelbox \\ [0.4cm]

\end{array}$
\end{center}
\caption{In the figure left, we see the minimal plane $\Sigma$ with $\PI \Sigma = \Gamma$. In the figure right, we see
the line segment $\beta$ which is transverse to $\Sigma$.}

\label{Figure 3}
\end{figure}

Now, we claim that there is another point $q \in \beta \cap \Sigma$ with $1<q_z<\delta$ which gives us a contradiction.
Let $N>0$ be such that $1+\frac{1}{N}<\delta$. Recall that $\Gamma_n = \Gamma_{n-1}\# C_n$ where $C_n$ is the round
circle in $\BR^2\times\{0\} \subset \Si$ with radius $1+\frac{1}{n}$ and center $(0,0,0)$. Let $\mu_1$ be the round
circle in $\BR^2\times\{0\} \subset \Si$ with radius $\frac{1}{N+1}<r_1<\frac{1}{N}$ and center $(0,0,0)$.  Let $\mu_2$
be the round circle in $\BR^2\times\{0\} \subset \Si$ with radius $\frac{1}{2N+1}<r_2<\frac{1}{2N}$ and center
$(0,0,0)$. Hence, $\mu_1$ is between the circles $C_N$ and $C_{N+1}$ in $\BR^2\times\{0\}$, and  $\mu_2$ is between the
circles $C_{2N}$ and $C_{2N+1}$ in $\BR^2\times\{0\}$. By choosing $r_1,r_2$ accordingly, further assume that $\mu_1
\cap \eta_N^\pm=\emptyset$ and $\mu_2\cap\eta_{2N}^\pm=\emptyset$ (See Figure 1). Let $P_i$ be the geodesic plane in
$\BH$ with $\PI P_i = \mu_i$, and let $\BH-P_i=\Omega^+_i\cup\Omega^-_i$ where $(0,0,0)\in\Omega^-_i$. Then, $P_1$ and
$P_2$ are least area planes in $X$, too.

Now, consider $\Sigma\cap P_i$. First, $\PI \Sigma=\Gamma$ and $\PI P_i = \mu_i$. $\mu_i$ intersect $\Gamma$ at exactly
2 points by construction. Say $\mu_i\cap \Gamma = \{x_i^+,x_i^-\}$. Since $\Sigma$ and $P_i$ are least area planes,
$\Sigma \cap P_i$ cannot contain a simple closed curve by Meeks-Yau exchange roundoff trick \cite{Co2}. Hence,
$\Sigma\cap P_i = \{l^i_j\}$ where $l^i_j$ is an infinite line segment in $\BH$ with $\PI l^i_j= \{x_i^+,x_i^-\}$.
Since $P_i$ is separating in $\BH$, all lines are separating in $\Sigma$. Hence, there is a natural ordering among
$\{l^i_j\}$. Let $l_1$ be the lowermost line among $\{l^1_j\}$ and let $l_2$ is the uppermost line among $\{l^2_j\}$
such that $l_1\cup l_2$ separates a component $\widehat{\Sigma}$ in $\Sigma$ where $\widehat{\Sigma}$ contains no line
segments $\{l^i_j\}$.

Now, $P_1$ and $P_2$ are the geodesic planes with $\PI P_i=\mu_i$, then let $\BH-P_i=\Omega^+_i\cup\Omega^-_i$ where
$(0,0,0)\in\Omega^-_i$. Let $Y= \overline{\Omega^+_1}\cap\overline{\Omega_2^-}$. By assumption on $r_1$ and $r_2$, we
know that $Y\cap \T_N =\emptyset$ and $Y\cap \T_{2N} =\emptyset$. Hence, $\widehat{Y} = Y\cap X$ would be a genus $N-1$
handlebody, i.e. $\widehat{Y} = Y- \bigcup_{i=N+1}^{2N-1} \T_i$. By construction, $\widehat\Sigma \subset \widehat{Y}$.
Topologically, we have a closed disk $\widehat{\Sigma}$ with $\partial \widehat{\Sigma} = \widehat{\alpha} \subset
\partial \widehat{Y}$. Hence, $\widehat{\alpha}$ is trivial element in $\pi_1(\widehat{Y})$. However, if we define
$\widehat{Z} = \widehat{Y} - \beta$, as in the proof of Lemma 2.4, we see that $\widehat{\alpha}$ is not a trivial
element in $\pi_1(\widehat{Z})$. This proves that $\widehat{\Sigma} \cap \beta \neq \emptyset$. Let $q$ be a point in
$\widehat{\Sigma} \cap \beta$. By construction $\frac{1}{2N}<q_z<\frac{1}{N}<\delta$. However, this contradicts with
the definition of $\delta$ as $\delta=\inf_{p\in\Sigma\cap\beta} \{\ p_z \in \BR \ | \ p=(p_x,p_y,p_z) \ \}$. The proof
follows.
\end{pf}

\section{Final Remarks}

We should note that our construction differs from the Freedman and He's heuristic construction in the following way. In
their construction, they want to apply to bridge principle to construct the sequence of minimal disks, then take the
limit. However, the examples of curves in $\Si$ described in \cite{La} shows that such a limit might not give a
connected plane, and bridges might escape to infinity. In our construction, it can be thought that we are still using
the bridges, but we are also using the tunnels acting as barrier which prevents bridges to escape to infinity. However,
because of these tunnels, while $\Sigma$ is a least area plane in $X$ ($\BH$ with tunnels deleted), it is not a least
area plane in $\BH$ anymore by \cite{Co1}.

On the other hand, one might try to use a "bridge principle at infinity" to construct such an example. In other words,
one might start with infinite family of geodesic planes as in this paper, and try to build bridges in $\Si$ to connect
the asymptotic boundaries of the geodesic planes, and take the limit. However, the problem with this approach would be
when you make a bridge at infinity, one might completely lost the original geodesic planes which goes through the
compact part, and get a completely different least area plane with the new asymptotic boundary which stays close to the
asymptotic sphere $\Si$. Hence, the barrier tunnels in our construction are very essential to construct such an
example.

Note also that it is known that if $\Sigma$ is a least area plane in $\BH$ with $\PI \Sigma = \Gamma$ where $\Gamma$ is
a simple closed curve which contains at least one smooth point, then $\Sigma$ is properly embedded in $\BH$ by
\cite{Co1}. However, since neither $\Sigma$ is not least area in our example nor $\Gamma$ is a simple closed curve,
\cite{Co1} does not apply here. Also, even though our example gives a complete, non-properly embedded minimal plane in
$\BH$, it is still not known the existence of a complete, non-properly embedded least area plane in $\BH$. So, it would
be an interesting question whether there exists a non-properly embedded least area plane in $\BH$.

As it is stated in the introduction, the key lemma of \cite{MR1} to prove its main result that {\em a complete,
embedded minimal surfaces with finite genus in $\BR^3$ is proper} is also true for $3$-manifolds with non-positive
curvature. However, the example constructed in this paper shows that even though the key lemma (Theorem 1 in
\cite{MR1}) is valid for minimal surfaces in $\BH$, it does not imply the properly embeddedness in $\BH$ like in
$\BR^3$.

If one considers our example constructed in this paper in the \cite{MR1} context, we get the following picture.
$\Sigma$ (in Theorem 2.5) is a minimal plane in $\BH$ with positive injectivity radius. If one applies Theorem 1 (or
Theorem 3) of \cite{MR1} to $\Sigma$ in $\BH$, we get  a lamination $\sigma = \overline{\Sigma} = \Sigma \cup P$ where
$P$ is the geodesic plane whose asymptotic boundary is the unit circle in $\BR^2\times\{0\}\subset \Si$. By using
Theorem 1, they prove Theorem 2 in \cite{MR1} which states that {\em a complete embedded connected minimal surface in
$\BR^3$ with positive injectivity radius is always properly embedded}.

Our example $\Sigma$ shows that Theorem 2 of \cite{MR1} is not true in $\BH$, whereas Theorem 1 of \cite{MR1} is valid
in $\BH$. Now, what goes wrong to get properly embeddedness of $\Sigma$ in $\BH$ as in the case of $\BR^3$ in
\cite{MR1}? In the proof of Theorem 2 in \cite{MR1}, Meeks and Rosenberg apply Theorem 1 to a complete embedded minimal
surface with positive injectivity radius in $\BR^3$, and in the closure, they get a minimal lamination $\mathcal{L}$.
By \cite{MR2}, they concluded that the limit leaves must be planes in $\BR^3$. Similarly, the limit leaf $P$ in our
lamination $\sigma$ is also a plane in $\BH$. So, everything is similar so far. However, when you apply Theorem 4 to
$\mathcal{L}$, they show that $M$ must have bounded curvature in an $\epsilon$ neighborhood of the limit leaf. However,
this contradicts to Lemma 1.3 of \cite{MR2} which states that $M$ cannot have unbounded curvature in the neighborhood
of a limit leaf. On the other hand, Theorem 4 is valid for $\sigma$ in $\BH$, too. Hence, we get that $\Sigma$ must
have bounded curvature in $\epsilon$ neighborhood of the limit leaf $P$. Unlike $\BR^3$, this can happen in $\BH$ case
as the high curvature regions which corresponds to {\em the bridges} in $\Sigma$ are far away from the limit leaf in
$\BH$. So, an analogous result of Lemma 1.3 in \cite{MR2} is not true in $\BH$ in general, and this is the place where
the technique in \cite{MR1} breaks down in $\BH$ case.

Note also that in \cite{MT}, Meeks and Tinaglia recently announced examples of non-properly embedded constant mean
curvature surfaces of finite topology for any $H\in [0,1)$ in $\BH$ . They also show that if $H\geq 1$ then the surface
must be properly embedded in $\BH$. Their example is different than ours, as they construct an infinite strip which is
a constant mean curvature surface limiting into two constant mean curvature annuli in $\BH$. The asymptotic boundary of
this surface is a pair of infinite lines where each line spirals into a pair of circles (asymptotic boundaries of the
annuli).

\end{document}